\documentclass[11pt]{article}

\usepackage{latexsym,amsfonts,amsmath,amssymb,amsthm}
\usepackage[english]{babel}
\usepackage[latin1]{inputenc}
\usepackage[T1]{fontenc}
\usepackage{graphics}
\usepackage{color}
\usepackage{url}

\newtheorem{duge}{Lemma}[section]
\newtheorem{prop}[duge]{Proposition}
\newtheorem{defi}[duge]{Definition}
\newtheorem{theo}[duge]{Theorem}

\newtheorem{cor}[duge]{Corollary}

\newcommand{\mcp}{\mathbb{P}}
\newcommand{\mce}{\mathbb{E}}

\newcommand{\mcn}{\mathbb{N}}
\newcommand{\mcr}{\mathbb{R}}
\newcommand{\un}{{\mathchoice {\rm 1\mskip-4mu l} {\rm 1\mskip-4mu l} {\rm 1\mskip-4.5mu l} {\rm 1\mskip-5mu l}}}

\newcommand{\srlun}{\mathcal{S}_1}
\newcommand{\srl}{\mathcal{S}}
\newcommand{\sr}{\mathcal{S^{\downarrow}}}

\newcommand{\pb}{\mathbb{P}^{(B)}}
\newcommand{\px}{\mathbb{P}^{(X)}}
\newcommand{\ex}{\mathbb{E}^{(X)}}
\newcommand{\eb}{\mathbb{E}^{(B)}}

\begin{document}

\title{On the equivalence of some eternal additive
coalescents}

\author{
Anne-Laure Basdevant}
\date{}
\maketitle

\begin{center}
\it{Laboratoire de Probabilités et Modèles Aléatoires,\\
 Universit\'e Pierre et Marie Curie,\\ 175 rue du Chevaleret,
75013 Paris, France.}
\end{center}

\vspace*{0.8cm}

\begin{abstract}
In this paper, we study additive coalescents. Using their
representation as fragmentation processes, we
 prove that the law of a large class of eternal additive
coalescents is absolutely continuous with respect to the law of the
standard additive coalescent on any bounded time interval.
\end{abstract}


\bigskip
\noindent{\bf A.M.S. Classification. } 60 J 25, 60 G 09.

\bigskip
\noindent{\bf e-mail. } Anne-Laure.Basdevant@ens.fr

\vspace*{1cm}

\section{Introduction}
The paper deals with additive coalescent processes, a class of
Markov processes which have been  introduced first by Evans and
Pitman \cite{Evanspitman98}. In the simple situation of a system
initially composed of a finite number $k$ of clusters with masses
$m_1,m_2,\ldots,m_k$,  the dynamics are such that each pair of
clusters $(m_i,m_j)$ merges into a unique cluster with mass
$m_i+m_j$ at rate $m_i+m_j$, independently of the other pairs. In
the sequel, we
 always assume that we start with a total mass equal to 1 (i.e.
$m_1+\ldots+m_k=1$). This induces no loss of generality  since we
can then deduce the law of any additive coalescent process through a
 time renormalization. Hence, an additive coalescent lives on the compact set
$$\sr=\{x=(x_i)_{i\ge 1}, x_1\ge x_2\ge \ldots \ge 0, \sum_{i\ge 1} x_i\le
1\},$$
 endowed with the topology of uniform
convergence.

 Evans and Pitman \cite{Evanspitman98} proved that we can
define an additive coalescent on the whole real line for a system
starting at time $t=-\infty$ with an infinite number of
infinitesimally small clusters. Such a process will be called an
eternal coalescent process. More precisely, if we denote by
$(C^n(t),t\ge 0)$ the additive coalescent starting from the
configuration $(1/n,1/n,\ldots,1/n)$, they proved that the sequence
of processes $(C^n(t+\frac{1}{2}\ln n),t\ge -\frac{1}{2}\ln n )$
converges  in distribution on the space of càdlàg paths with values
in the set $\sr$
 toward some process
$(C^\infty(t),t \in \mcr)$, which  is called the standard additive
coalescent. We stress that this process is defined for all time
$t\in \mcr$. A remarkable property of the standard additive
coalescent is that, up to time-reversal, its becomes a fragmentation
process. Namely, the process $(F(t),t\ge 0)$ defined by
$F(t)=C^{\infty}(-\ln t)$ is a self-similar fragmentation process
with index of self similarity $\alpha=1/2$, with no erosion and with
dislocation measure $\nu$
given by
$$\nu(x_1\in dy)=(2\pi y^3(1-y)^3)^{-1/2}dy \quad \mbox{for } y\in]1/2,1[,\quad \nu(x_3>0)=0.$$
We refer to Bertoin \cite{CoursBertoin03} for the definition of
erosion, dislocation measure, and index of self similarity of a
fragmentation process and a proof. Just recall that in a
fragmentation process, distinct fragments evolve independently of
each others.

 Aldous and Pitman \cite{Aldouspitman98}  constructed
this fragmentation process $(F(t),t\ge 0)$  by cutting the skeleton
of the continuum Brownian random tree according to a Poisson point
process. In another paper \cite{Aldouspitman00}, they gave a
generalization of this result: consider for each $n\in \mcn$ a
decreasing sequence $r_{n,1}\ge \ldots\ge r_{n,n}\ge 0$ with sum 1,
set $\sigma^2_n=\sum_{i=1}^n r_{n,i}^2$ and suppose that
$$\lim_{n\rightarrow \infty}{\sigma_n}=0 \mbox{ and }
\lim_{n\rightarrow\infty} \frac{r_{n,i}}{\sigma_n}=\theta_i \mbox{
for all }i\in \mcn.$$ Assume further that $\sum_i{\theta_i^2}<1$ or
$\sum_i{\theta_i}=\infty$. Then, it is proved in
\cite{Aldouspitman00}
  that if $M^n=(M^n(t),t\ge 0)$ denotes the additive coalescent
process starting with $n$ clusters with mass $r_{n,1}\ge \ldots\ge
r_{n,n}$, then $(M^{(n)}(t-\ln \sigma_n),t\ge\ln\sigma_n)$ has a
limit distribution as $n\rightarrow \infty$, which can be obtained
by cutting a specific inhomogeneous random tree with a point Poisson
process. Furthermore, any extreme eternal additive coalescent can be
obtained  this way up to a deterministic time translation.

Bertoin \cite{Bertoin01b} gave another construction of the limit of
the process  $(M^{(n)}(t-\ln \sigma_n),t\ge\ln\sigma_n)$ in the
following way. Let $b_\theta$ be the bridge with exchangeable
increments defined for $s\in [0,1]$ by
$$b_\theta(s)=\sigma b_s+\sum_{i=1}^{\infty}\theta_i(\un_{\{s\ge
V_i\}}-s),$$ where $(b_s,s\in[0,1])$ is a standard Brownian bridge,
$(V_i)_{i\ge 1}$ is an i.i.d. sequence of uniform random variable on
[0,1] independent of $b$ and $\sigma=1-\sum_i \theta_i^2$. Let
 $\varepsilon_\theta=(\varepsilon_\theta(s),s\in [0,1])$ be the excursion obtained from $b_\theta$ by Vervaat's transform, i.e.
$\varepsilon_\theta(s)=b_{\theta}(s+m \mod 1)-b_{\theta}(m),$ where
$m$ is the point of [0,1] where $b_{\theta}$ reaches its minimum.
For all $t\ge 0$, consider
$$\varepsilon^{(t)}_\theta(s)=ts-\varepsilon_\theta(s), \hspace*{1cm} S^{(t)}_\theta(s)=\sup_{0\le u\le s
}\varepsilon^{(t)}_\theta(u),$$ and define $F^\theta(t)$ as the
sequence of the lengths of the constancy intervals of the process
$(S^{(t)}_\theta(s),0\le s\le 1)$. Then  the limit of the process
$(M^{(n)}(t-\ln \sigma_n),t\ge\ln\sigma_n)$ has the law of
$(F^\theta(e^{-t}),t\in \mcr)$.
 Miermont \cite{Miermont01} studied the same process in the special
 case where $\varepsilon_\theta$ is the normalized excursion above the minimum of a
 spectrally negative Lévy process. More precisely  let
$(X_t,t\ge 0)$ be a Lévy process with no positive jump, with
unbounded variation and with positive and finite  mean. Let
$\overline{X}(t)=\sup_{0\le s\le t} X_t$ and denote by
$\varepsilon_X=(\varepsilon_X(s),s\in [0,1])$  the normalized
excursion with duration 1 of the reflected process $\overline{X}-X$.
We now define in the same way as for $b_\theta$, the processes
$\varepsilon_X^{(t)}(s)$, $S^{(t)}_X(s)$ and $F^X(t)$. Then, the
process $(F^X(e^{-t}),t\in \mcr)$ is a mixture of some eternal
additive coalescents (see \cite{Miermont01} for more details).
Furthermore, $(F^X(t),t\ge 0)$ is a fragmentation process in the
sense that  distinct fragments evolve independently of each other
(however, it is not necessarily homogeneous in time). It is quite
remarkable that the Lévy property of $X$ ensures the branching
property of $F^X$. We stress that there exist other eternal additive
coalescents for which this property fails. Notice that when  the
Lévy process $X$ is the standard Brownian motion $B$, the process
$(F^{B}(e^{-t}),t\in \mcr)$ is then the standard additive coalescent
and $(F^{B}(t),t\ge 0)$ is a self-similar and time-homogeneous
fragmentation process.

In this paper, we study the relationship  between the laws $\px$  of
$(F^{X}(t),t\ge 0)$ and $\pb$ of  $(F^{B}(t),t\ge 0)$.  We prove
that, for certain Lévy processes $(X_t,t\ge 0)$, the law $\px$ is
absolutely continuous with respect to $\pb$ and we compute
explicitly the density. Our main result is the following:

\begin{theo}\label{absolu}
Let $(\Gamma(t),t\ge 0)$ be a subordinator with no drift. Assume
that $\mce(\Gamma_1)<\infty$ and take any $c\ge \mce(\Gamma_1)$. We
define  $X_t=B_t-\Gamma_t+ct$, where $B$ denotes a Brownian motion
independent of $\Gamma$. Let $(p_t(u),u\in \mcr)$ and $(q_t(u),u\in
\mcr)$ stand for the respective density of $B_t$ and $X_t$. In
particular $p_t(u)=\frac{1}{\sqrt{2\pi t}}\exp(-\frac{u^2}{2t})$.
Let $\srlun$ be the space of positive sequences with sum 1. We
consider the function $\mathbf{h}:\mcr_+\times \srlun$ defined by
$$\mathbf{h}(t,\mathbf{x})=e^{tc}\frac{p_1(0)}{q_1(0)}\prod_{i=1}^{\infty}\frac{q_{x_i}(-tx_i)}{p_{x_i}(-tx_i)}
 \qquad \mbox{ with } \;\mathbf{x}=(x_i)_{i\ge 1}.$$
Then, for all $t\ge 0$, the function $\mathbf{h}(t,\cdot)$ is
bounded on $\srlun$ and has the following properties:
\begin{itemize}
\item $\mathbf{h}(t,F(t))$ is a $\pb$-martingale,
\item for every $t\ge 0$, the law of the process $(F^X(s),0\le s\le t)$ is absolutely continuous with respect
to that of $(F^B(s),0\le s\le t)$ with density
$\mathbf{h}(t,F^B(t)).$
\end{itemize}
\end{theo}

Let us  notice that $\mathbf{h}(t,\cdot)$ is a multiplicative
function, i.e. it can be written as the product of functions, each
of them depending only on the size of a single fragment.  In the
sequel we will use the notation
$$h(t,x)=
e^{tcx}\left(\frac{p_1(0)}{q_1(0)}\right)^x\frac{q_{x}(-tx)}{p_{x}(-tx)}
 \qquad \mbox{ for }x\in ]0,1] \mbox{ and } t\ge 0,$$
 so we have $\mathbf{h}(t,\mathbf{x})=\prod_i h(t,x_i)$.  This multiplicative form  of
$\mathbf{h}(t,\cdot)$ implies that the process $F^X$ has the
branching property (i.e. distinct fragments evolve independently of
each other) since $F^B$ has it. Indeed, for every multiplicative
bounded continuous function $\mathbf{f}:\sr\mapsto \mcr_+$, for all
$t'>t>0$ and $\mathbf{x}\in \sr$, we have, since
$\mathbf{h}(t,F^B(t))$ is a $\pb$-martingale,
$$\ex\Big(\mathbf{f}(F(t'))\,\big|\,F(t)=\mathbf{x}\Big)=\frac{1}{\mathbf{h}(t,\mathbf{x})}\eb\Big(\mathbf{h}(t',F(t'))\mathbf{f}(F(t'))\,\big|\,F(t)=\mathbf{x}\Big).$$
Using the branching property of $F^B$ and the multiplicative form of
$\mathbf{h}(t,\cdot)$,  we get
$$\ex\Big(\mathbf{f}(F(t'))\,\big|\,F(t)=\mathbf{x}\Big)=\frac{1}{\mathbf{h}(t,\mathbf{x})}\prod_{i}\eb\Big(\mathbf{h}(t',F(t'))\mathbf{f}(F(t'))\,\big|\,F(t)=(x_i,0,\ldots)\Big).$$
And finally we deduce
\begin{eqnarray*}\ex\Big(\mathbf{f}(F(t'))\,\big|\,F(t)=\mathbf{x}\Big)&=&\frac{1}{\mathbf{h}(t,\mathbf{x})}\prod_{i}h(t,x_i)
\ex\Big(\mathbf{f}(F(t'))\,\big|\,F(t)=(x_i,0,\ldots)\Big)\\
&=&\prod_{i}\ex\Big(\mathbf{f}(F(t'))\,\big|\,F(t)=(x_i,0,\ldots)\Big).
\end{eqnarray*}
Let $M_{\mathbf{x}}$ (resp. $M_{x_i}$) be the random measure on
]0,1[ defined by $M_{\mathbf{x}}=\sum_i \delta_{s_i}$ where the
sequence $(s_i)_{i\ge 1}$ has the law of $F(t')$ conditioned on
$F(t)=\mathbf{x}$ (resp. $F(t)=(x_i,0,\ldots)$). Hence we have, for
every bounded continuous function $g:\mcr \mapsto \mcr$,
$$\mce\Big(\exp(-<g,M_{\mathbf{x}}>)\Big)=\prod_{i=1}^{\infty}\mce\Big(\exp(-<g,M_{x_i}>)\Big),$$
which proves that $M_{\mathbf{x}}$ has the law of $\sum_i M_{x_i}$
where the random measures $(M_{x_i})_{i\ge 1}$ are independent.
Hence the process $F^{X}$ has the branching property. Notice also
that other multiplicative martingales have already been studied in
the case of branching random walks
\cite{Biggins77,Chauvinrouault88,Neveu88,Kyprianou04}.

\vspace*{0.5cm}

This paper will be divided in two sections. The first section is
devoted to the proof of this theorem and in the next one, we will
use the fact that $\mathbf{h}(t,F^B(t))$ is a $\pb$-martingale to
describe an integro-differential equation solved by the function
$h$.

\section{Proof of Theorem \ref{absolu}}
The assumptions and notation in Theorem \ref{absolu} are implicitly
enforced throughout this section.
\subsection{Absolute continuity}

In order to prove Theorem \ref{absolu}, we will first prove the
absolute continuity of the law  $\px_t$ of $F^{X}(t)$  with respect
to the law  $\pb_t$ of  $F^{X}(t)$ for a fixed time $t>0$ and for a
finite number of fragments. We begin first by a definition:

\begin{defi} Let
$x=(x_1,x_2,\ldots)$ be a sequence of positive numbers with sum 1.
We call the random variable $y=(x_{j_1},x_{{j_2}},\ldots)$ a size
biased rearrangement of $x$ if we have:
$$\forall i\in \mcn, \;\mcp(j_1=i)= x_i,$$ and by induction
$$\forall i\in \mcn\backslash\{i_1,\ldots,i_k\},\; \mcp(j_{k+1}=i \; |\;
j_1={i_1},\ldots,j_k={i_k})=\frac{x_i}{1-\sum_{l=1}^k x_{i_l}}.$$

\end{defi}

Notice that  for every Lévy process $X$ satisfying hypotheses of
Theorem \ref{absolu}, we have $\sum_{i=1}^{\infty}F_i(t)=1 \quad
\px_t\mbox{-a.s.}$ (it is clear by the construction  from an
excursion of $X$ since $X$ has unbounded variation, cf
\cite{Miermont01}, Section 3.2). Hence the above definition  can be
applied to $F^{X}(t)$.

The following lemma gives the distribution of the first  $n$
fragments of $F^{X}(t)$, chosen with a size-biased pick:
\begin{duge}
Let $(\tilde{F}^X_1(t),\tilde{F}^X_2(t),\ldots)$ be a size biased
rearrangement of $F^X(t)$. Then for all $n\in \mcn$, for all
$x_1,\ldots,x_n \in \mcr_+$  such that $S=\sum_{i=1}^n x_i < 1$, we
have
$$ \px_t(\tilde{F}^X_1\in dx_1,\ldots, \tilde{F}^X_n\in
dx_n)=\frac{t^n}{q_1(0)}q_{1-S}(St)\prod_{i=1}^n
\frac{q_{x_i}(-tx_i)}{1-\sum_{k=1}^i x_k}dx_1\ldots dx_n.$$
\end{duge}

 \begin{proof} On the one hand, Miermont \cite{Miermont01} gave a description of the law of $F^X(t)$:
 let $T^{(t)}$ be a subordinator with Lévy measure
  $z^{-1}q_z(-tz)\un_{z>0}dz$.
 Then $F^X(t)$ has the law of the sequence of the jumps of $T^{(t)}$ before time $t$
 conditioned on $T^{(t)}_t=1$.

 One the other hand, consider a subordinator $T$ on the time interval $[0,u]$   conditioned by
 $T_u=y$ and pick a jump of $T$ by size-biased sampling. Then, its
 distribution has  density
$$\frac{zuh(z)f_u(y-z)}{yf_u(y)}dz,$$
where $h$ is the density of the Lévy measure of $T$ and $f_u$ is the
density of $T_u$ (see Theorem 2.1 of \cite{Permanpityor92}). Then,
in the present case, we have
$$u=t,\quad y=1,\quad h(z)=z^{-1}q_z(-tz), \quad f_u(z)=\frac{u}{z}q_z(u-zt) \quad \mbox{ (cf. Lemma 9 of  \cite{Miermont01})}.$$
Hence we get $$\px_t(\tilde{F}^X_1\in
dz)=\frac{tq_z(-tz)q_{1-z}(zt)}{(1-z)q_1(0)}dz.$$ This proves the
lemma in the case $n=1$. The proof for $n\ge 2$ uses an induction.
Assume  that we have proved the case $n-1$ and let us prove the case
$n$. We have
\begin{multline*}
$$\px_t(\tilde{F}^X_1\in dx_1,\ldots, \tilde{F}^X_n\in
dx_n)=\\
\px_t(\tilde{F}^X_1\in dx_1,\ldots, \tilde{F}^X_{n-1}\in
dx_{n-1})\px_t(\tilde{F}^X_n\in dx_n\;|\;\tilde{F}^X_1\in
dx_1,\ldots, \tilde{F}^X_{n-1}\in dx_{n-1}).$$
\end{multline*}
Furthermore, Perman, Pitman and Yor \cite{Permanpityor92} have
proved that the $n$-th size biased picked jump $\Delta_n$ of a
subordinator before time $u$ conditioned by $T_u=y$ and
$\Delta_1=x_1,\ldots, \Delta_{n-1}=x_{n-1}$ has the law of a size
biased picked jump  of the subordinator $T$ before time $u$
conditioned by $T_u=y-x_1-\ldots - x_{n-1}.$ Hence we get:
\begin{multline*}
$$\px_t(\tilde{F}^X_1\in dx_1,\ldots, \tilde{F}^X_n\in
dx_n)=\\
\left(\frac{t^{n-1}}{q_1(0)}q_{1-S_{n-1}}(S_{n-1}t)\prod_{i=1}^{n-1}
\frac{q_{x_i}(-tx_i)}{1-S_i}\right)
\frac{tq_{x_n}(-tx_n)q_{1-S_{n}}(S_{n}t)}{(1-S_n)q_{1-S_{n-1}}(S_{n-1}t)}dx_1\ldots
dx_{n},$$
\end{multline*}
 where $S_i=\sum_{k=1}^i x_k$. And so the lemma is proved by induction.

 \end{proof}

Since the lemma is clearly also true for $\pb$ (take $\Gamma=c=0$),
we get:
\begin{cor}\label{rapportfini} Let $(F(t),t\ge 0)$ be a fragmentation process. Let
$(\tilde{F}_1(t),\tilde{F}_2(t),\ldots)$ be a size biased
rearrangement of $F(t)$. Then for all $n\in \mcn$, for all
$x_1,\ldots,x_n \in \mcr_+$  such that $S=\sum_{i=1}^n x_i < 1$, we
have
$$\frac{\px_t(\tilde{F}_1\in dx_1,\ldots, \tilde{F}_n\in dx_n)}{\pb_t(\tilde{F}_1\in
dx_1,\ldots, \tilde{F}_n\in dx_n)}=h_n(t,x_1,\ldots,x_n),$$
$$ \mbox{with } h_n(t,x_1,\ldots,x_n)=
\frac{p_1(0)}{q_1(0)}\frac{q_{1-S}(St)}{p_{1-S}(St)}\prod_{i=1}^n
\frac{q_{x_i}(-tx_i)}{p_{x_i}(-tx_i)}.
$$

\end{cor}

To establish that  the law of $F^X(t)$ is  absolutely  continuous
 with respect to the law of $F^B(t)$ with density $\mathbf{h}(t,\cdot)$, it remains to check that the
function $h_n$ converges as $n$ tends to infinity to $\mathbf{h}$
$\pb_t$-a.s. and in $L^1(\pb_t)$. In this direction, we first prove
two lemmas:

\begin{duge}\label{produit} We have
$\frac{q_{y}(-ty)}{p_{y}(-ty)}<1$ for all $y>0$ sufficiently small.
As a consequence, if $(x_i)_{i\ge 1}$ is a sequence of positive
numbers with $\lim_{i\rightarrow \infty} x_i= 0$, then the product
$\prod_{i=1}^{n}\frac{q_{x_i}(-tx_i)}{p_{x_i}(-tx_i)}$ converges as
$n$ tends to infinity.
\end{duge}

\begin{proof}
Since $X_t=B_t-\Gamma_t+tc$, notice that we have
$$\forall s>0,\; \forall u \in \mcr,\quad q_s(u)=\mce\Big(p_s(u+\Gamma_s-cs)\Big).$$
Hence if we replace  $p_s(u)$ by its expression $\frac{1}{\sqrt{2\pi
s}}\exp(-\frac{u^2}{2s})$, we get
\begin{equation}\label{qu}
\frac{q_s(u)}{p_s(u)}=\exp\left(cu-\frac{c^2s}{2}\right)\mce\left[\exp\left(-\frac{\Gamma^2_s}{2s}-\Gamma_s(\frac{u}{s}-c)\right)\right].
\end{equation}
 i.e., for all $y>
0$, for all $t\ge 0$,
$$\frac{q_{y}(-ty)}{p_{y}(-ty)}=\exp\left(-y(ct+\frac{c^2}{2})\right)\mce\left[\exp\left(-\frac{\Gamma^2_y}{2y}+\Gamma_y(t+c)\right)\right].$$
Using  the inequality $(c-a)(c-b)\ge -\left(\frac{b-a}{2}\right)^2$,
we have
$$-\frac{\Gamma^2_y}{2y}+\Gamma_y(t+c)\le \frac{y(t+c)^2}{2}$$ and we deduce $$\frac{q_{y}(-ty)}{p_{y}(-ty)}\le
e^{\frac{t^2 y}{2}}.$$

Fix $c'\in]0,c[$, let $f$ be the function defined by
$f(y)=\mcp(\Gamma_y\le c'y)$. Since $\Gamma_t$ is a subordinator
with no drift, we have $\lim_{y\rightarrow 0}{f(y)}=1$ (indeed,
$\Gamma_y=o(y)$ a.s., see \cite{Bertoin96}). On the event
$\{\Gamma_y\le c'y\}$, we have
\begin{eqnarray*}\exp\left(-y(ct+\frac{c^2}{2})\right)\exp\left(-\frac{\Gamma^2_y}{2y}+\Gamma_y(t+c)\right)&\!\le\!&
\exp(-y(\frac{1}{2}(c-c')^2+t(c-c')))\\&\!\le\!& \exp(-\varepsilon
y),\end{eqnarray*}
 with $\varepsilon=\frac{1}{2}(c-c')^2.$ Hence, we
get the upper bound
$$\frac{q_{y}(-ty)}{p_{y}(-ty)}\le e^{-\varepsilon
y}f(y)+(1-f(y))e^{\frac{yt^2}{2}}.$$ Since $f(y)\rightarrow 1$ as
$y\rightarrow 0$, we deduce
$$e^{-\varepsilon
y}f(y)+(1-f(y))e^{\frac{yt^2}{2}}=1-\varepsilon y+ o(y).$$ Thus, we
have $\frac{q_{y}(-ty)}{p_{y}(-ty)}<1$ for $y$ small enough, and so
the product converges for every sequence $(x_i)_{i\ge 0}$ which
tends to $0$.
\end{proof}

We prove now a second lemma:
\begin{duge}\label{S} We have
$$\lim_{s\rightarrow 1^-} \frac{q_{1-s}{(st)}}{p_{1-s}(st)}=e^{tc}.$$
\end{duge}

\begin{proof}
We use again Identity (\ref{qu}) established in the proof of Lemma
\ref{produit}. We get:
$$\frac{q_{1-s}(st)}{p_{1-s}(st)}=\exp\left(tsc-\frac{c^2}{2}(1-s)\right)
\mce\left[\exp\left(-\frac{\Gamma_{1-s}^2}{2(1-s)}-\Gamma_{1-s}(\frac{ts}{1-s}-c)\right)\right].$$
For $s$ close enough to 1, $\frac{ts}{1-s}-c\ge 0$, hence we get
$$\mce\left[\exp\left(-\frac{\Gamma_{1-s}^2}{2(1-s)}-\Gamma_{1-s}(\frac{ts}{1-s}-c)\right)\right]\le
1$$ and we deduce $$\limsup_{s\rightarrow 1^-}
\frac{q_{1-s}(st)}{p_{1-s}(st)}\le e^{tc}.$$

For the lower bound, we write
\begin{multline*}
\mce\!\left[\exp\!\left(-\frac{\Gamma_{1-s}^2}{2(1-s)}-\Gamma_{1-s}(\frac{ts}{1-s}-c)\right)\!\right]\\
\begin{aligned}
&\ge
\mce\!\left[\exp\!\left(-\frac{\Gamma_{1-s}}{2(1-s)}-\Gamma_{1-s}(\frac{ts}{1-s}-c)\!\right)\un_{\{\Gamma_{1-s}\le
1\}}\right]\\
&\ge\mce\!\left[\exp\left(-\Gamma_{1-s}\frac{1+2ts}{2(1-s)}\right)\un_{\{\Gamma_{1-s}\le
1\}}\right]\\
&\ge\mce\!\left[\exp\left(-\Gamma_{1-s}\frac{1+2ts}{2(1-s)}\right)\right]-\mcp(\Gamma_{1-s}\ge
1).\\\end{aligned}
\end{multline*}
Since $\Gamma_t$ is a subordinator with no drift,
$\lim_{u\rightarrow 0}\frac{\Gamma_u}{u}=0$ a.s., and we have for
all $K>0$,
$$\lim_{u\rightarrow 0^+}\mce\left[\exp\left(-K\frac{\Gamma_u}{u}\right)\right]=1.$$
Hence, we get
$$\liminf_{s\rightarrow 1^-} \frac{q_{1-s}(st)}{p_{1-s}(st)}\ge e^{tc}.$$
\end{proof}

We are now able to prove the absolute continuity of  $\px_t$ with
respect to  $\pb_t$. Since $S_n=\sum_{i=1}^{n}x_i$ converges
$\pb_t$-a.s. to 1, Lemma \ref{produit} and \ref{S} imply that
$H_n=h_n(t,\tilde{F}_1(t),\ldots,\tilde{F}_n(t))$ converges to
$H=\mathbf{h}(t,F(t))$ $\pb$-a.s.

Let us now prove that $H_n$ is uniformly bounded, which implies the
$L^1$ convergence. We have already proved that there exists
$\varepsilon>0$ such that:
$$\forall x\in ]0,\varepsilon[,\quad \frac{q_{x}(-tx)}{p_{x}(-tx)}\le 1.$$
Besides, it is well known that, if $X_t=B_t-\Gamma_t+ct$, its
density $(t,u)\rightarrow q_t(u)$ is continuous on $\mcr_+^*\times
\mcr$. Hence, on $[\varepsilon,1]$, the function $x\rightarrow
\frac{q_{x}(-tx)}{p_{x}(-tx)}$ is continuous and we can find an
upper bound $A>0$ of this function . As there are at most
$\frac{1}{\varepsilon}$ fragments of $F(t)$ larger than
$\varepsilon$, we deduce the upper bound:
$$\prod_{i=1}^\infty
\frac{q_{F_i}(-tF_i)}{p_{F_i}(-tF_i)}\le
A^{\frac{1}{\varepsilon}}.$$ Likewise, the function $S\rightarrow
\frac{q_{1-S}(St)}{p_{1-S}(St)}$ is continuous on $[0,1[$ and has a
limit at $1$, so it is bounded by some $D>0$ on $[0,1]$. Hence we
get
$$H_n\le A^{\frac{1}{\varepsilon}} D \frac{p_1(0)}{q_1(0)} \quad\pb \mbox{-a.s.}$$
So $H_n$ converges to $H\quad\pb \mbox{-a.s.}$ and in $L^1(\pb)$.
Furthermore, by construction, $H_n$ is a $\pb$-martingale, hence we
get for all $n\in\mcn$,
$$\mce^{(B)}(H\;|\;\tilde{F}_1,\ldots,\tilde{F}_n)=H_n,$$
and so, for every bounded continuous function $f:\srlun\rightarrow
\mcr$, we have
$$\mce^{(X)}\Big[f(F(t))\Big]=\mce^{(B)}\Big[f(F(t))\mathbf{h}(t,F(t))\Big].$$

Hence, we have proved that, for a fixed time $t\ge 0$, the law of
$F^X(t)$ is absolutely continuous with respect to that of $F^B(t)$
with density $\mathbf{h}(t,F^B(t))$. Furthermore, Miermont
\cite{Miermont01}  has proved that the processes $(F^X(e^{-t}),t\in
\mcr)$ and $(F^B(e^{-t}),t\in \mcr)$ are both eternal additive
coalescents (with different entrance laws). Hence, they have the
same semi-group of transition and  we get the absolute continuity of
the law of the process $(F^X(s),0\le s\le t)$ with  respect to that
of $(F^B(s),0\le s\le t)$ with density $\mathbf{h}(t,F^B(t))$.

\subsection{Sufficient condition for equivalence} We can
now wonder whether the measure $\px$ is equivalent to  the measure
$\pb$, that is whether $\mathbf{h}(t,F(t))$ is strictly positive
$\pb$-a.s. A sufficient condition is given by the following
proposition.

\begin{prop}\label{proppositif} Let $\phi$ be the  Laplace exponent of the subordinator $\Gamma$,
i.e. $$ \forall s\ge 0, \forall q\ge 0,\quad
\mce(\exp(-q\Gamma_s))=\exp(-s\phi(q)).$$ Assume that there exists
$\delta>0$ such that
\begin{equation}\label{positif}\lim_{x\rightarrow \infty}
\phi(x)x^{\delta-1}=0,
\end{equation} then the function $\mathbf{h}(t,F(t))$ defined
in Theorem \ref{absolu} is strictly positive $\pb$-a.s.
\end{prop}

We stress that the condition \ref{positif} is very weak. For
instance, let $\pi$ be the Lévy measure of the subordinator and
$I(x)=\int_{0}^x \overline{\pi}(t)dt$ where $\overline{\pi}(t)$
denotes $\pi(]t,\infty[)$. It is well known that $\phi(x)$ behaves
like $xI(1/x)$ as $x$ tends to infinity (see \cite{Bertoin96}
Section III). Thus, the condition \ref{positif} is equivalent to
$I(x)=o(x^\delta)$ as $x$ tends to 0 (recall that we always have
$I(x)=o(1)$).

\begin{proof}
Let $t>0$.  We must check that
$\prod_{i=1}^{\infty}\frac{q_{x_i}(-tx_i)}{p_{x_i}(-tx_i)}$ is
$\pb_t$-almost surely strictly positive. Using (\ref{qu}), we have:

$$\frac{q_{y}(-ty)}{p_{y}(-ty)}=\exp\left(-y(ct+\frac{c^2}{2})\right)\mce\left[\exp\left(-\frac{\Gamma^2_y}{2y}+\Gamma_y(t+c)\right)\right].$$
Since we have $\sum_{i=1}^{\infty} x_i=1 \quad \pb_t$-a.s., we get
$$\prod_{i=1}^{\infty}\frac{q_{x_i }(-tx_i)}{p_{x_i}(-tx_i)}\ge \exp\left(-ct+\frac{c^2}{2}\right)
\prod_{i=1}^{\infty}\mce\left[\exp\left(-\frac{\Gamma^2_{x_i}}{2x_i}+c
\Gamma_{x_i} \right)\right].
$$
Hence we have to find a lower bound for
$\mce\left[\exp\left(-\frac{\Gamma^2_y}{2y}+c \Gamma_y
\right)\right]$. Since $c\ge \mce(\Gamma_1)$, we have
$$\mce\left[\exp\left(-\frac{\Gamma^2_y}{2y}+c \Gamma_y
\right)\right]\ge
\mce\left[\exp\left(\frac{\Gamma_y}{y}(\mce(\Gamma_y)-
\frac{\Gamma_y}{2}) \right)\right].$$
  Set  $A=\mce(\Gamma_1)$ and let us fix $K>0$. Notice that the event $\mce(\Gamma_y)-
\frac{\Gamma_y}{2}\ge -Ky$ is equivalent to  the event $\Gamma_y\le
(2A+K)y$ and by Markov inequality, we have
$$\mcp(\Gamma_y\ge (2A+K)y)\le \frac{A}{2A+K}.$$

Hence we get
\begin{eqnarray*}
\mce\left[\exp\left(-\frac{\Gamma^2_y}{2y}+c \Gamma_y
\right)\right]& \ge &
\mce\left[\exp\left(\frac{\Gamma_y}{y}(\mce(\Gamma_y)-
\frac{\Gamma_y}{2})\un_{\{\Gamma_y\le (2A+K)y\}} \right)\right]\\
&\ge & \mce\left(\exp(-K\Gamma_y)\un_{\{\Gamma_y\le (2A+K)y\}} \right)\\
&\ge & \mce\left(\exp(-K\Gamma_y)
\right)-\mce\left(\exp(-K\Gamma_y)\un_{\{\Gamma_y> (2A+K)y\}}
\right)\\
&\ge & \exp(-\phi(K)y)-\frac{A}{2A+K}.
\end{eqnarray*}
This inequality holds for all $K>0$.  Hence, with $\varepsilon>0$
and $K=y^{-\frac{1}{2}-\varepsilon}$, we get
$$\mce\left[\exp\left(-\frac{\Gamma^2_y}{2y}+c \Gamma_y
\right)\right] \ge
\exp\left(-\phi(y^{-\frac{1}{2}-\varepsilon})y\right)- A
y^{\frac{1}{2}+\varepsilon}.$$

Furthermore,  the product
$\prod_{i=1}^{\infty}\mce\left[\exp\left(-\frac{\Gamma^2_{x_i}}{2x_i}+c
 \Gamma_{x_i}
\right)\right]$ is strictly positive if  the series
$$\sum_{i=1}^\infty
1-\mce\left[\exp\left(-\frac{\Gamma^2_{x_i}}{2x_i}+c \Gamma_{x_i}
\right)\right]$$ converges. Hence, a sufficient condition is
$$\exists \;\varepsilon >0 \mbox{ such that } \sum_{i=1}^\infty\left(
1-\exp\left(-\phi(x_i^{-\frac{1}{2}-\varepsilon})x_i\right)+
x_i^{\frac{1}{2}+\varepsilon}\right)<\infty \quad
\pb_t\mbox{-a.s.}$$

Recall that the distribution of the Brownian fragmentation at time
$t$ is equal to the distribution of the  jumps of a stable
subordinator $T$ with index $1/2$ before time $t$ conditioned on
$T_t=1$ (see \cite{Aldouspitman98}). Hence, it is well known that we
have for all $\varepsilon>0$
$$ \sum_{i=1}^\infty x_i^{\frac{1}{2}+\varepsilon}<\infty \quad \pb_t\mbox{-a.s.} \quad\mbox{(see Formula (9) of \cite{Aldouspitman98})}.$$

Thus, we have equivalence between $\pb_t$ and $\px_t$ as soon as
there exist two strictly positive numbers $\varepsilon,\varepsilon'$
such that, for $x$ small enough
$$\phi(x^{-\frac{1}{2}-\varepsilon})x\le
x^{\frac{1}{2}+\varepsilon'}.$$ One can easily check that this
condition is equivalent to (\ref{positif}).
\end{proof}

\vspace*{0.5cm}
 In Theorem \ref{absolu}, we have supposed that
$X_t$ can be written as $B_t+\Gamma_t-ct$, with $c\ge
\mce(\Gamma_1)$ and $\Gamma_t$ subordinator. We can wonder whether
the theorem applies for a larger class of Lévy processes. Notice
first that the process $X$ must fulfill the conditions of Miermont's
paper \cite{Miermont01} recalled in the introduction, i.e. $X$  has
no positive jumps, unbounded variation and finite and positive mean.
Hence, a possible extension of the Theorem would be for example for
$X_t=\sigma^2B_t+\Gamma_t-ct$, with $\sigma>0$, $\sigma\neq 1$. In
fact, it is clear that Theorem \ref{absolu} fails in this case. Let
just consider for example $X_t=2B_t$. Using Proposition 3 of
\cite{Miermont01}, we get that
$$(F^X(2t),t\ge 0)\overset{\tiny{\mbox{law}}}{=}(F^B(t),t\ge 0).$$
But, it is well known that we have
$$\lim_{n\rightarrow \infty}n^2F^{\downarrow}_n(t)=t\sqrt{2/\pi} \quad
\pb\mbox{-a.s.} \quad \mbox{ (see \cite{Bertoin04}})
$$
Hence, the laws $\pb_t$ and $\pb_{2t}$ are mutually singular.

\section{An integro-differential equation}

Since $\mathbf{h}(t,F(t))$ is the density of $\px$ with respect to
$\pb$ on the sigma-field $\mathcal{F}_t=\sigma(F(s),s\le t)$, it is
a $\pb$-martingale. Hence, in this section, we will compute the
infinitesimal generator of a fragmentation to deduce a remarkable
integro-differential equation.

\subsection{The infinitesimal generator of a fragmentation process}
In this section, we recall a result obtained by Bertoin and Rouault
in an unpublished paper \cite{Bertoinrouault05}.

We denote by $\mathcal{D}$ the space of functions $f: [0,1]\mapsto
]0,1]$ of class $\mathcal{C}^1$ and with $f(0)=1$. For $f\in
\mathcal{D}$ and $\mathbf{x}\in \sr$, we set
$$\mathbf{f}(\mathbf{x})=\prod_{i=1}^{\infty}f(x_i).$$
For $\alpha\in \mcr_+$ and $\nu$ measure on $\sr$ such that
$\int_{\sr}(1-x_1)\nu(d\mathbf{x})<\infty$, we define the operator
$$G_{\alpha}\mathbf{f}(\mathbf{x})=\mathbf{f}(\mathbf{x})\sum_{i=1}^{\infty}x_i^\alpha\int\nu(d\mathbf{y})
\left(\frac{\mathbf{f}(x_i\mathbf{y})}{f(x_i)}-1\right)\quad \mbox{
for } f\in \mathcal{D} \mbox{ and } \mathbf{x}\in \sr.$$
\begin{prop}\label{generateur}
Let $(X(t),t\ge 0)$ be a self-similar fragmentation with index of
self-similarity $\alpha>0$, dislocation measure $\nu$ and no
erosion. Then, for every function $f\in \mathcal{D}$, the process
$$\mathbf{f}(X(t))-\int_{0}^{t}G_{\alpha}\mathbf{f}(X(s))ds$$ is a martingale.
\end{prop}

\begin{proof}
We will first prove the following lemma
\begin{duge}\label{lemmegen}
For $f \in \mathcal{D}, \mathbf{y}\in \sr, r\in[0,1]$, we have
$$\Big| \frac{\mathbf{f}(r\mathbf{y})}{f(r)}-1\Big|\le 2C_f e^{C_f} r(1-y_1),$$
with $C_f=\left|\left|\frac{f'}{f^{2}}\right|\right|_{\infty}$.
\end{duge}

Notice that, since $f$ is $\mathcal{C}^1$ on $[0,1]$ and strictly
positive, $C_f$ is always finite.
\begin{proof}
First, we write
$$|\ln f(ry_1)-\ln f(r)|\le \Big|\Big|\frac{f'}{f}\Big|\Big|_{\infty}(1-y_1)r\le
C_f(1-y_1)r.$$ We deduce then
 $$\frac{\mathbf{f}(r\mathbf{y})}{f(r)}-1\le
 \frac{f(ry_1)}{f(r)}-1\le e^{C_f(1-y_1)r}-1\le C_fe^{C_f}(1-y_1)r. $$
 Besides we have $$\ln \frac{1}{f(x_1)}\le \frac{1}{f(x_i)}-1\le
 C_f x_i,\quad
\mbox{ which implies }\quad\mathbf{f}(\mathbf{x})\ge
f(x_1)\exp(-C_f\sum_{i=2}^{\infty}x_i).$$ Hence we get
$$\frac{\mathbf{f}(r\mathbf{y})}{f(r)}\ge
\frac{f(ry_1)}{f(r)}\exp(-C_f(1-y_1)r)\ge \exp(-2C_f(1-y_1)r),$$ and
we deduce $$1-\frac{\mathbf{f}(r\mathbf{y})}{f(r)}\le 2C_f e^{C_f}
(1-y_1)r. $$
\end{proof}
We can now prove Proposition \ref{generateur}.
 We denote by
$\mathcal{T}$ the set of times where some dislocation occurs (which
is a countable set). Hence we can write
$$\mathbf{f}(X(t))-\mathbf{f}(X(0))=\sum_{s\in [0,t]\cap \mathcal{T}}\Big(
\mathbf{f}(X(s))-\mathbf{f}(X(s-))\Big),$$ as soon as
$$\sum_{s\in [0,t]\cap \mathcal{T}}\Big|
\mathbf{f}(X(s))-\mathbf{f}(X(s-))\Big|<\infty$$
  For $s\in \mathcal{T}$,
if the $i$-th fragment $X_i(s-)$ is involved in the dislocation, we
set $k_s=i$ and we denote by $\Delta_s$ the element of $\sr$
according to $X(s-)$ has been broken. Hence, we have
$$\sum_{s\in [0,t]\cap \mathcal{T}}\Big|
\mathbf{f}(X(s))-\mathbf{f}(X(s-))\Big|=\sum_{ s\in \mathcal{T}\cap
[0,t]}\mathbf{f}(X(s-))\left(\sum_{i=1}^\infty
\un_{k_s=i}\Big|\frac{\mathbf{f}(X_i(s-)\Delta_s)}{f(X_i(s-))}-1\Big|\right).$$
Hence, since a fragment of mass $r$ has a rate of dislocation
$\nu_r(dx)=r^{\alpha}\nu(dx)$,  the predictable compensator is
\begin{multline*}\int_{0}^t ds \;\mathbf{f}(X(s-))\int_{\sr}
\nu(d\mathbf{y})\sum_{i=1}^\infty X^\alpha_i(s-)
\Big|\frac{\mathbf{f}(X_i(s-)\mathbf{y})}{f(X_i(s-))}-1\Big|\\
\begin{aligned}
&\le 2C_fe^{C_f}\int_{0}^t  \sum_{i=1}^\infty X_i(s-)
\int_{\sr}(1-y_1)\nu(d\mathbf{y})ds.\\
&\le 2C_fe^{C_f}t \int_{\sr}(1-y_1)\nu(d\mathbf{y})
\end{aligned}
\end{multline*}

Hence $$\sum_{s\in [0,t]\cap \mathcal{T}}\Big|
\mathbf{f}(X(s))-\mathbf{f}(X(s-))\Big|<\infty \quad \mbox{a.s.},$$
and thus we have
$$\mathbf{f}(X(t))-\mathbf{f}(X(0))=\sum_{s\in [0,t]\cap \mathcal{T}}\Big(
\mathbf{f}(X(s))-\mathbf{f}(X(s-))\Big),$$ i.e.
$$\mathbf{f}(X(t))-\mathbf{f}(X(0))=\sum_{ s\in \mathcal{T}\cap
[0,t]}\mathbf{f}(X(s-))\left(\sum_{i=1}^\infty
\un_{k_s=i}\left(\frac{\mathbf{f}(X_i(s-)\Delta_s)}{f(X_i(s-))}-1\right)\right),$$
whose predictable compensator is
$$\int_{0}^t ds \;\mathbf{f}(X(s-))\int_{\sr}
\nu(d\mathbf{y})\sum_{i=1}^\infty X^\alpha_i(s-)
\left(\frac{\mathbf{f}(X_i(s-)\mathbf{y})}{f(X_i(s-))}-1\right)=\int_0^t
G_{\alpha}\mathbf{f}(X(s))ds.$$
\end{proof}

\subsection{Application to $\mathbf{h}(t,F(t))$}
Let $F(t)$ be a fragmentation process and $q_t(x)$ be the density of
a Lévy process fulfilling the hypotheses of Theorem \ref{absolu}. We
have proved in the first section that the function
$$H_t=\mathbf{h}(t,F(t))=e^{tc}\frac{p_1(0)}{q_1(0)}\prod_{i=1}^{\infty}\frac{q_{F_i(t)}(-tF_i(t))}{p_{F_i(t)}(-tF_i(t))}$$
is a $\pb$-martingale (since it is equal to
$\frac{d\px}{d\pb}|\mathcal{F}_t$). We set
$$g(t,x)=e^{tcx}\frac{q_{x}(-tx)}{p_{x}(-tx)} \quad \mbox{for } x\in]0,1], t\ge 0\quad \mbox { and } g(t,0)=1.$$
$$\mbox{ Set now }\quad
\mathbf{g}(t,\mathbf{x})=\prod_{i=1}^{\infty}g(t,x_i(t)) \quad
\mbox{for }\mathbf{x}\in \sr,t\ge 0.$$ So we have, as $\sum_i
F_i(t)=1$ $\pb$-a.s.,
$$H_t=\frac{p_1(0)}{q_1(0)}\mathbf{g}(t,F(t)) \quad \mbox{ for all } t\ge 0.$$

It is well known that if $q_t(u)$ is the density of a Lévy process
$X_t=B_t-\Gamma_t+ct$, the function $(t,u)\mapsto q_t(u)$ is
$\mathcal{C}^{\infty}$ on $\mcr_+^*\times \mcr$. Hence $(t,x)\mapsto
g(t,x)$ is also $\mathcal{C}^{\infty}$ on $\mcr_+\times ]0,1]$ and
in particular,  for all $x\in [0,1]$, the function $t\rightarrow
g(t,x)$ is $\mathcal{C}^1$ and so $\partial_t g(t,x)$ is well
defined. The next proposition gives a integro-differential equation
solved by the function $g$ when $g$ has some properties of
regularity at points $(t,0)$, $t\in \mcr_+$.

\begin{prop}\label{propmartingale}
\begin{enumerate}
\item \label{point1}
Assume that  for all $t\ge 0$, $\partial_x g(t,0)$ exists and the
function
 $(t,x)\rightarrow
\partial_x g(t,x)$ is continuous at $(t,0)$. Then $g$
solves the equation:
$$\left\{\begin{array}{l}\partial_t
 g(t,x)+ \sqrt{x}\displaystyle\int_{0}^1\frac{dy}{\sqrt{8\pi
 y^3(1-y)^3}}\Big(g(t,xy)g(t,x(1-y))-g(t,x)\Big)=0\\
 g(0,x)=\frac{q_{x}(0)}{p_{x}(0)}.\end{array}\right.$$
\item \label{point2} If the Lévy measure of the subordinator $\Gamma$ is finite, then the
above conditions  on $g$  hold.
\end{enumerate}
\end{prop}

\begin{proof}Let us first notice  that the hypotheses of the proposition
imply that the integral
$$\displaystyle\int_{0}^1\frac{dy}{\sqrt{8\pi
 y^3(1-y)^3}}\Big(g(t,xy)g(t,x(1-y))-g(t,x)\Big)$$
is well defined and is continuous in $x$ and in $t$. Indeed, this
integral is equal to
$$2\displaystyle\int_{0}^{\frac{1}{2}}\frac{dy}{\sqrt{8\pi
 y^3(1-y)^3}}\Big(g(t,xy)g(t,x(1-y))-g(t,x)\Big).$$
And for all $y\in ]0,1/2[, x\in]0,1], t\in \mcr_+$, there exist
$c,c'\in [0,x]$ such that
$$\frac{g(t,xy)g(t,x(1-y))-g(t,x)}{y}=x(g(t,x)\partial_x g(t,c)-g(t,xy)\partial_x
g(t,c')).$$ Thanks to the hypothesis  that the function
$(t,x)\rightarrow
\partial_x g(t,x)$ is continuous on $\mcr_+\times [0,1]$, $|x(g(t,x)\partial_x g(t,c)-g(t,xy)\partial_x
g(t,c'))|$ is uniformly bounded on $[0,T]\times [0,1]\times
[0,\frac{1}{2}]$ and so by application of the theorem of dominated
convergence, the integral is continuous in $t$ on $\mcr_+$ and in
$x$ on [0,1].

 We begin by proving the first point  of the
proposition. Recall that, according to Proposition \ref{generateur},
the generator of the Brownian fragmentation is
$$G_{\frac{1}{2}}\mathbf{f}(\mathbf{x})=\mathbf{f}(\mathbf{x})\sum_{i=1}^{\infty}\sqrt{x_i}\int\nu(d\mathbf{y})
\left(\frac{\mathbf{f}(x_i\mathbf{y})}{f(x_i)}-1\right) ,$$ with
$$\nu(y_1\in du)=(2\pi u^3(1-u)^3)^{-1/2}du \quad \mbox{for } u\in]1/2,1[,\quad \nu(y_1+y_2\neq 1)=0 \quad\mbox{(cf. \cite{Bertoin02})}.$$
  Hence,
$$M_t=\mathbf{g}(t,F(t))-\mathbf{g}(0,F(0))-\int_{0}^{t}G_{\frac{1}{2}}\mathbf{g}(s,F(s))+\partial_t
\mathbf{g}(s,F(s))ds$$ is a $\pb$-martingale. Since
$\mathbf{g}(t,F(t))$ is already a $\pb$-martingale, we get
$$G_{\frac{1}{2}} \mathbf{g}(s,F(s))+\partial_t \mathbf{g}(s,F(s))=0  \quad\pb\mbox{-a.s.}\quad\mbox{ for almost every }
s>0,$$ i.e. for almost every  $s>0$

 $$\mathbf{g}(s,F(s))\sum_{i=1}^{\infty}\left[
 F^{1/2}_i(s)\int_{\sr}\nu(d\mathbf{y})\left(\frac{\mathbf{g}(s,F_i(s)\mathbf{y})}{g(s,F_i(s))}-1\right)+\frac{\partial_t g(s,F_i(s))}{g(s,F_i(s))}\right]=0
  \quad\pb\mbox{-a.s.}$$
With $F(s)=(x_1,x_2,\ldots)$, we get
$$\sum_{i=1}^{\infty}\left[
 x^{1/2}_i\int_{\sr}\nu(d\mathbf{y})\left(\frac{\mathbf{g}(s,x_i\mathbf{y})}{g(s,x_i)}-1\right)+\frac{\partial_t
 g(s,x_i)}{g(s,x_i)}\right]=0 \quad \pb_s\mbox{-a.s.}$$
Notice also that this series is absolutely convergent. Indeed,
thanks to Lemma \ref{lemmegen}, we have
$$\Big|x^{1/2}_i\int_{\sr}\nu(d\mathbf{y})\left(\frac{\mathbf{g}(s,x_i\mathbf{y})}{g(s,x_i)}-1\right)\Big|\le C_{g,s}x_i \int_{\sr}(1-y_1)\nu(d\mathbf{y}),$$
where $C_{g,s}$ is a positive constant (which depends on $g$ and
$s$), and, besides we have
$$g(t,x)=\exp\left(-x\frac{c^2}{2}\right)\mce\left[\exp\left(-\frac{\Gamma^2_x}{2x}+\Gamma_x(t+c)\right)\right].$$
Thus, by application of the theorem of dominated convergence, it is
easy to prove that the function $t\rightarrow
\mce\left[\exp\left(-\frac{\Gamma^2_x}{2x}+\Gamma_x(t+c)\right)\right]$
is derivable with derivative
$$\partial_t
\mce\left[\exp\left(-\frac{\Gamma^2_x}{2x}+\Gamma_x(t+c)\right)\right]=
\mce\left[\Gamma_x\exp\left(-\frac{\Gamma^2_x}{2x}+\Gamma_x(t+c)\right)\right].$$
Notice also that this quantity is continuous in $x$ on [0,1].

  Hence
we have
 $$\forall x_i\in]0,1[, \forall s>0, \quad\frac{\partial_t
 g(s,x_i)}{g(s,x_i)}>0.$$
Thus we deduce $$\sum_{i=1}^\infty \frac{\partial_t
 g(s,x_i)}{g(s,x_i)}<\infty \quad \pb_s\mbox{-a.s.}$$
 Let define $$k(t,x)=\partial_t g(t,x)+ \sqrt{x}\int_{0}^1\frac{dy}{\sqrt{8\pi
 y^3(1-y)^3}}\Big(g(t,xy)g(t,x(1-y))-g(t,x)\Big).$$
 Hence we have
 \begin{equation}\sum_{i=1}^{\infty} k(s,x_i)=0 \quad\pb_s\mbox{-a.s.}\quad\mbox{ for almost every }
s>0,\end{equation} and
\begin{equation}\sum_{i=1}^{\infty} |k(s,x_i)|<\infty \quad\pb_s\mbox{-a.s.}\quad\mbox{ for almost every }
s>0.\end{equation} Furthermore, $x\rightarrow k(t,x)$ is continuous
on $[0,1]$, hence, thanks to the following lemma, we get for almost
every $s>0$, $k(s,x)=0$ for $x\in [0,1]$.  And, since $s\rightarrow
k(s,x)$ is continuous on $\mcr_+$, we deduce $k\equiv 0$ on
$\mcr_+\times [0,1]$.

\end{proof}

\begin{duge}
Fix $t>0$. Let $\pb_t$ denote the law of the Brownian fragmentation
at time $t$. Let $k : [0,1]\mapsto \mcr$ be a continuous  function,
  such that
$$\sum_{i=1}^{\infty} k(x_i)=0 \quad \pb_t\mbox{-a.s. and } \sum_{i=1}^{\infty}| k(x_i)|<\infty \quad \pb_t\mbox{-a.s.}$$
Then $k\equiv 0$ on [0,1].
\end{duge}

\begin{proof}
Let $F(t)=(F_1(t),F_2(t)\ldots)$ be a Brownian fragmentation at time
$t$ where the sequence $(F_i(t))_{i\ge 1}$ is ordered by a
size-biased pick. We denote by $\srl$ the set of positive sequence
with sum less than 1. Since $F(t)$ has the law of the size biased
reordering of the jumps of a stable subordinator $T$ (with index
$1/2$) before time $t$, conditioned by $T_t=1$ (see
\cite{Aldouspitman98}), it is obvious that we have
$$\forall x\in
]0,1-S[,\quad \pb_t(F_1\in dx\;| \;(F_i)_{i\ge 3})>0,  $$ where
$S=\sum_{i\ge 3} F_i$. Let $\mathbb{Q}_t$ be the measure on $\srl$
defined by
$$\forall A\subset \srl, \quad \mathbb{Q}_t(A)=\pb_t((F_i)_{i\ge 3}\in
A)$$ and $\lambda$ the Lebesgue measure on $[0,1]$. Hence we have,
for all $y\in \srl$ - $\mathbb{Q}_t$-a.s.
$$ \forall x\in ]0,S[,\quad
k(x)+k(1-S-x)+\sum_{i=1}^{\infty}k(y_i)=0
\quad\lambda\mbox{-a.s.},$$ where $S=\sum_i y_i$. We choose now
$y\in \srl$ such that this equality holds for almost every $x\in
]0,S[$. Thus, we get that there exists a constant $C=C(y)$ such that
$$k(x)+k(1-S-x)=C,\quad \mbox{ for all } x\in ]0,S[\quad
\lambda\mbox{-a.s.}$$ Since $k$ is continuous, this equality holds
in fact for all $x\in [0,S]$. Furthermore, we have also
$$\forall s\in ]0,1[, \quad \mathbb{Q}_t(S\in ds)>0.$$
Hence, this implies the existence  for almost every $s\in ]0,1[$ of
a constant $C_s$ such that
$$k(x)+k(1-s-x)=C_s\quad \mbox{ for all } x\in ]0,s[.$$
Thanks to the continuity of $k$, we can deduce that this property
holds in fact for all $s\in [0,1]$. Hence we have
$$\forall x,y\in [0,1]^2, \mbox{ such that } x+y\le 1,\;
k(x+y)=k(x)+k(y).$$ So $k$ is a linear function and since
$\sum_{i=1}^{\infty} x_i=1$  $\;\pb_t\mbox{-a.s.}$, we get $k\equiv
0$ on [0,1].
\end{proof}

We prove now the point \ref{point2} of Proposition
\ref{propmartingale}.

\begin{proof} Assume that the Lévy measure of $\Gamma$ is finite.
It is obvious that $g$ has the same regularity that the function
$\frac{q_{x}(-tx)}{p_{x}(-tx)}$. Recall now that we have
$$\frac{q_{x}(-tx)}{p_{x}(-tx)}=\exp\left(-x(ct+\frac{c^2}{2})\right)\mce\left[\exp\left(-\frac{\Gamma^2_x}{2x}+\Gamma_x(t+c)\right)\right].$$
 Hence a
sufficient condition for $g$ to fulfill the hypotheses of
Proposition \ref{propmartingale} is
\begin{itemize}
\item $u_t(x)=
\mce\left[\exp\left(-\frac{\Gamma^2_x}{2x}+\Gamma_x(t+c)\right)\right]$
is derivable at 0.
\item $w(t,x)=u_t'(x)$ is continuous at $(t,0)$ for $t\in \mcr_+$.
\end{itemize}
We write $u_t(x)=a_t(x,x)$ with
$$a_t(y,z)=\mce\left[\exp\left(-\frac{\Gamma^2_y}{2z}+\Gamma_y(t+c)\right)\right].$$
Since  the function $(y,z)\rightarrow
\frac{y^2}{2z^2}\exp\left(-\frac{y^2}{2z}+y(t+c)\right)$ is bounded
on $\mcr_+\times[0,1]$, we get
$$\partial_z
a_t(y,z)=\mce\left[\frac{\Gamma^2_y}{2z^2}\exp\left(-\frac{\Gamma^2_y}{2z}+\Gamma_y(t+c)\right)\right]
\quad \mbox{ for } z\in ]0,1].$$ Recall that the generator of a
subordinator with no drift and Lévy measure $\pi$ is given  for
every  bounded function $f$ $C^1$ with bounded derivative  by
$$ \forall y\in \mcr_+,\;
Lf(y)=\int_{0}^{\infty}\!(f(y+s)-f(y))\pi(ds), \quad (\mbox{c.f.
Section 31 of \cite{Sato99}}).$$ Hence, we get for all $z_0>0$,
\begin{multline*}
$$\partial_y
a_t(y,z_0)=\mce(L a_t(\Gamma_y,z_0))\\=\mce\left[\int_{0}^{\infty}
\left(\exp\left(-\frac{(\Gamma_y+s)^2}{2z_0}+(\Gamma_y+s)(t+c)\right)-
\exp\left(-\frac{\Gamma_y^2}{2z_0}+\Gamma_y(t+c)\right)\right)\pi(ds)\right].$$
\end{multline*}
And we deduce
\begin{multline*}
$$u_t'(x)=\mce\left[\frac{\Gamma^2_x}{2x^2}\exp\left(-\frac{\Gamma^2_x}{2x}+\Gamma_x(t+c)\right)\right]\\+\mce\left[\int_{0}^{\infty}
\left(\exp\left(-\frac{(\Gamma_x+y)^2}{2x}+(\Gamma_x+y)(t+c)\right)-
\exp\left(-\frac{\Gamma_x^2}{2x}+\Gamma_x(t+c)\right)\right)\pi(dy)\right],
$$
\end{multline*}

We must prove that $(t,x)\rightarrow u_t'(x)$ is continuous at
$(t,0)$ for $t\ge 0$.
 For every Lévy measure $\pi$, the first term has
limit 0 as $(t',x)$ tends to $(t,0)$ (by dominated convergence). For
the second term, notice that we  have for all $x\in ]0,1]$,
$$\Big|\exp\left(\!-\frac{(\Gamma_x+y)^2}{2x}+(\Gamma_x+y)(t+c)\!\right)-
\exp\left(\!-\frac{\Gamma_x^2}{2x}+\Gamma_x(t+c)\!\right)\Big|\le
2\exp\left(\!\frac{(t+c)^2x}{2}\!\right),$$ and for all $y>0$,
$\exp\left(\!-\frac{(\Gamma_x+y)^2}{2x}+(\Gamma_x+y)(t+c)\!\right)-
\exp\left(\!-\frac{\Gamma_x^2}{2x}+\Gamma_x(t+c)\!\right)$ converges
almost surely to $-1$ as $(t',x)$ tends to $(t,0)$. Hence,  if
$\pi(\mcr_+)<\infty$, we deduce that the $\lim_{(t',x)\rightarrow
(t,0)} u_t'(x)$  exists (and is equal to $-\pi(\mcr_+)$).
\end{proof}
\nocite{Bertoin96}
\bibliographystyle{plain}

\bibliography{biblioall}
\end{document}